\documentclass[11pt]{article}
\usepackage{verbatim} 
\usepackage{amsmath}
\usepackage{amssymb}
\usepackage{amsfonts}
\usepackage{amsthm}
\usepackage{cite}
\usepackage{graphicx} 
\usepackage{pgf,tikz}
\usepackage{tikz-cd} 
\usepackage{mathrsfs}
\usetikzlibrary{arrows}
\usepackage{color}
\usepackage[legalpaper, margin=1in]{geometry}
\usepackage{enumitem}
\usepackage{latexsym}

\newtheorem{theorem}{Theorem}[section]
\newtheorem*{theorem*}{Theorem}

\newtheorem{corollary}[theorem]{Corollary}
\newtheorem{lemma}[theorem]{Lemma}
\newtheorem{proposition}[theorem]{Proposition}

\theoremstyle{definition}

\newtheorem{question}[theorem]{Question}
\newtheorem{remark}[theorem]{Remark}
\newtheorem{claim}{Claim}[theorem]
\theoremstyle{remark}
\newtheorem*{claimproof}{Proof}

\newcommand{\R}{\mathbb{R}}
\newcommand{\N}{\mathbb{N}}

\newcommand{\mesh}{\mathrm{mesh}}
\newcommand{\claimend}{{\hfill $\blacksquare$}}

\renewcommand{\restriction}{\mathord{\upharpoonright}}

\title{Compact spaces homeomorphic to their respective squares} 
\author{Jan Dud\'ak\footnote{https://orcid.org/0000-0003-0627-6641}
\footnote{Supported by the grant GACR 24-10705S}\\
Department of Mathematical Analysis\\ 
Faculty of Mathematics and Physics, Charles University\\
Prague, Czechia\\
E-mail: dudakjan@seznam.cz
\and 
Benjamin Vejnar\footnote{https://orcid.org/0000-0002-2833-5385} \footnote{
Supported by the grant GACR 24-10705S
}\\
Department of Mathematical Analysis\\ 
Faculty of Mathematics and Physics, Charles University\\
Prague, Czechia\\
E-mail: vejnar@karlin.mff.cuni.cz}

\begin{document}
\maketitle

\renewcommand{\thefootnote}{}



\renewcommand{\thefootnote}{\arabic{footnote}}
\setcounter{footnote}{0}

To the memory of W. J. Charatonik.

\begin{abstract}
We deal with topological spaces homeomorphic to their respective squares. Primarily, we investigate the existence of large families of such spaces in some subclasses of compact metrizable spaces. As our main result we show that there is a family of size continuum of pairwise non-homeomorphic compact metrizable zero-dimensional spaces homeomorphic to their respective squares. This answers a question of W. J. Charatonik. We also discuss the situation in the classes of continua, Peano continua and absolute retracts.
\end{abstract}

\section{Introduction}

There are many (separable metrizable) topological spaces which are known to be homeomorphic to their respective squares. Examples include the rationals, the irrationals, the Cantor space, the Erdös complete or rational space, an infinite discrete space, an infinite power of any space, mutual products of these spaces, etc. In some of these cases, a topological characterization of the space allows us to verify that the space is homeomorphic to its square.
If we restrict ourselves to the compact metrizable case, the topological dimension becomes a natural obstruction since if $\dim X\geq1$ for a compact metrizable space $X$, then $\dim X^n\geq n$ for each natural number $n$ \cite[1.8.K (b)]{EngelkingDimension}. Hence, a compact metrizable space homeomorphic to its square is necessarily either zero-dimensional or infinite-dimensional.


On May 14 2020 W. J. Charatonik presented during the Wroclaw Set Theory seminar a joint result with S. Sahan about the existence of uncountably many pairwise non-homeomorphic zero-dimensional compact metrizable spaces homeomorphic to their respective squares \cite{CharatonikSahan}. In his presentation, he asked whether such a collection can be found of cardinality $\mathfrak{c}=|\mathbb R|$. In this paper we answer the question positively, that is, we prove the following theorem.

\begin{theorem}\label{MainResult}
    There exists a family $\mathcal{F}$ of size continuum of pairwise non-homeomorphic compact metrizable zero-dimensional spaces such that $X \times X$ is homeomorphic to $X$ for each $X \in \mathcal{F}$.
\end{theorem}

A naive attempt to prove this theorem by searching for spaces of the form $X=Y^{\mathbb N}$ is doomed to fail since if $Y^{\mathbb N}$ is a compact metrizable zero-dimensional space with more than one point, then it is homeomorphic to the Cantor space. It should also be noted that a countable compact space is homeomorphic to its square if and only if it has at most one point, since the Cantor-Bendixson rank of any infinite countable metrizable compact space $X$ is strictly less than that of $X \times X$.

Obviously, if a topological space $X$ is homeomorphic to $X \times X$, then it is homeomorphic to $X^n$ for each $n \in \N$. On the other hand, if $X$ is homeomorphic to $X^n$ for some $n>2$, it may not be the case that $X$ is homeomorphic to $X^2$. A class $\mathcal C$ of topological spaces is said to have the \emph{Tarski cube property} if every $X\in \mathcal C$ which is homeomorphic to its cube is homeomorphic to its square. Trnková proved in \cite{Trnkova1980CMUC} that the class of countable metrizable spaces has the Tarski cube property while the class of connected metrizable spaces does not. Going further, it was shown in \cite{Trnkova1980PAMS} that every compact zero-dimensional metrizable space $X$ homeomorphic to $X^n$ for some $n \in \N \setminus \{ 1 \}$ is actually homeomorphic to $X^2$. It is worth noting that the original statement of the last mentioned result was formulated in the language of Boolean algebras, their free products and isomorphisms, and can be translated via the Stone duality to the language of compact zero-dimensional spaces, their products and homeomorphisms.

It was shown by van Douwen in \cite[Theorem 17.1]{vanDouwen} that every locally compact space homeomorphic to its square has a compactification homeomorphic to its square. The fundamental idea behind the construction can be generalized into a universal method of obtaining spaces homeomorphic to their respective square as follows. Starting with any topological space $X$ and any continuous mapping $f:X^2\to X$ consider the inverse sequence
\[
\begin{tikzcd}[column sep = huge]
X & X^2 \arrow{l}[swap]{f} & (X^2)^2 \arrow{l}[swap]{f \times f} & (X^2)^4 \arrow{l}[swap]{f \times f \times f \times f} & \cdots \arrow[l]
\end{tikzcd}
.\]
Then the inverse limit is a space homeomorphic to its square. Conversely, every space $X$ homeomorphic to its square can be realized using this method in a trivial way (take any homeomorphism $f \colon X^2 \to X$, then the inverse limit of the above inverse sequence is a space homeomorphic to X). The Pe{\l}{}czy\'{n}ski compactum is a nice example of a space which can be easily obtained with this method in a non-trivial way: If $X=\N \cup \{ \infty \}$ is the one-point compactification of $\N$ and $f:X^2\to X$ maps $\N^2$ onto $\N$ injectively and $X^2 \setminus \N^2$ onto $\{ \infty \}$, then the inverse limit of the above inverse sequence is a space homeomorphic to the Pe{\l}{}czy\'{n}ski compactum.

There exists a Peano continuum which is homeomorphic to its square but not to its infinite powers \cite[p.87-89]{vanMill}. The proof relies, among other things, on the idea behind van Douwen's result.

Let us mention one more result somewhat related to our topic. Recently, Medini and Zdomskyy proved that every filter in $P(\omega)/fin$ (considered as a subspace of $P(\omega)\sim 2^{\omega}$) is homeomorphic to its square \cite{MediniZdomskyy}.

Let us briefly outline the content of this paper: In Section 2 we define preliminary notions. In Section 3 we discuss various ways to produce continuum many pairwise non-homeomorphic infinite dimensional compact metrizable spaces which are homeomorphic to their respective squares. We also discuss the cases of Peano continua, absolute retracts and some other classes. In Section 4 we prove Theorem \ref{MainResult} as the main result of this paper. In Section 5 we provide a correction of the proof of \cite[Theorem 3.3]{CharatonikSahan}, which is a theorem very important for us since its corollary, \cite[Theorem 3.4]{CharatonikSahan}, is used in the proof of our Theorem \ref{MainResult} in Section 4. 

\section{Preliminaries}

We use the symbol $\mathfrak{c}$ for the cardinality of the continuum, that is, $\mathfrak{c} = |\R |$. We denote by $\omega$ the least infinite ordinal and by $\omega_1$ the least uncountable ordinal. Ordinal numbers and ordinal arithmetic operations will play a role in section \ref{SectionZeroDim}.

For any topological space $X$ and any subset $S$ of $X$, we denote by $S'$ the set of all the points in $X$ which are limit points of $S$. A subset $P$ of a topological space is said to be perfect if it is closed and has no isolated points. Thus, $P$ is perfect if and only if $P'=P$.

\begin{lemma}\label{SubstractingOneIsolatedPointDoesntChangeAnything}
    Let $X$ be a metrizable space with infinitely many isolated points and let $z \in X \setminus X'$. Then $X \setminus \{ z \}$ is homeomorphic to $X$.
\end{lemma}

\begin{proof}
    If $X'$ is not open in $X$, fix an infinite set $A \subseteq X \setminus X'$ containing $z$ such that $A$ has exactly one limit point. If $X'$ is open in $X$, let $A:=X \setminus X'$. Either way, fix any bijection $\varphi \colon A \to A \setminus \{ z \}$ and define a mapping $h \colon X \to X \setminus \{ z \}$ by $h(x)=\varphi(x)$ for $x \in A$ and by $h(x)=x$ for $x \in X \setminus A$. It is easy to see that $h$ is a homeomorphism.
\end{proof}

A topological space is said to be Polish if it is separable and completely metrizable. A subspace $Y$ of a Polish space $X$ is Polish if and only if $Y$ is $G_\delta$ in $X$ (see, e.g., \cite[Theorem 3.11]{Kechris}).

By \cite[Theorem 6.2]{Kechris}, every nonempty perfect subset of a Polish space contains a homeomorphic copy of the Cantor space. In particular, every nonempty countable Polish space has an isolated point. By the Cantor-Bendixson theorem (see, e.g., \cite[Theorem 6.4]{Kechris}), for every Polish space $X$, there is a unique perfect set $P \subseteq X$ such that $X \setminus P$ is countable. We call this set the perfect kernel of $X$ and we denote it by $\mathrm{PK}(X)$. The complement of $\mathrm{PK}(X)$ can be characterized as the maximal countable open subset of $X$. In other words, $X \setminus \mathrm{PK}(X)$ is the set of all the points in $X$ which have a countable neighbourhood. Note that $\mathrm{PK}(X)$ is nonempty (and thus of size continuum) if and only if $X$ is uncountable. These facts and observations imply the following proposition.

\begin{proposition}\label{IsolatedPointsAreDenseInTheComplementOfThePerfectKernel}
     For every Polish space $X$, the set $X \setminus \mathrm{PK}(X)$ is contained in the closure of the set $X \setminus X'$. In particular, if $X$ is countable, then $X \setminus X'$ is dense in $X$.
\end{proposition}

For a metrizable compact space $X$, the Cantor-Bendixson derivative of order $\alpha$ of $X$, denoted by $X^{(\alpha)}$, is defined for every ordinal number $\alpha$ inductively as follows:
\begin{enumerate}[label=(\arabic*),noitemsep]
    \item $X^{(0)}=X$,
    \item $X^{(\alpha)}= \big( X^{(\beta)} \big)'$ if $\alpha = \beta +1$,
    \item $X^{(\alpha)}= \bigcap \big\lbrace X^{(\gamma)} \, ; \ \gamma < \alpha \big\rbrace$ if $\alpha$ is a limit ordinal.
\end{enumerate}
If $X$ is countable, there exists an ordinal $\alpha < \omega_1$ with $X^{(\alpha+1)}=\emptyset$. The least such $\alpha$ is called the Cantor-Bendixson rank of $X$ and we denote it by $\mathrm{CB}(X)$. If $X$ is, in addition, nonempty, $\mathrm{CB}(X)$ can be characterized as the unique ordinal $\alpha$ for which $X^{(\alpha)}$ is nonempty and finite.

For every metrizable compact space $X$ and every $x \in X \setminus \mathrm{PK}(X)$, we let
\[ \mathrm{CB}(x,X):=\min \big\lbrace \mathrm{CB}(K) \, ; \ K \subseteq X \textup{ is a countable compact neighbourhood of } x \big\rbrace. \]
It is easy to show that $\mathrm{CB}(x,X)$ is the least ordinal $\alpha$ for which $x \notin X^{(\alpha+1)}$.

For any metrizable topological space $X$, we denote by $\mathcal{K}(X)$ the hyperspace of nonempty compact subsets of $X$ and we endow it with the Vietoris topology. If $(X,d)$ is a metric space, we denote by $d_H$ the Hausdorff metric on $\mathcal{K}(X)$ induced by $d$. That is, for all $K,L \in \mathcal{K}(X)$,
\[ d_H (K,L) = \max \Big\lbrace \max_{x \in K} \mathrm{dist}_d (x,L) \, , \ \max_{y \in L} \mathrm{dist}_d (y,K) \Big\rbrace. \]
The topology induced by $d_H$ is exactly the Vietoris topology on $\mathcal{K}(X)$. Also, recall that a metrizable topological space $X$ is compact if and only if $\mathcal{K}(X)$ is compact.

For any topological space $X$ and subsets $Y$ and $Z$ of $X$ with $Z \subseteq Y$, a partition $\mathcal{A}$ of $Z$ is said to be $Y$-clopen provided that every set $A \in \mathcal{A}$ is relatively clopen in $Y$.

For any family $\mathcal{A}$ of subsets of a metric space $X$, we denote $\mathrm{mesh} (\mathcal{A}) := \sup \big\lbrace \mathrm{diam}(A) \, ; \ A \in \mathcal{A} \big\rbrace$.

\section{Infinite dimensional case} 

In this section we discuss several ways to construct a large family of (pairwise non-homeomorphic) continua homeomorphic to their respective squares.

In 1934 Waraszkiewicz constructed an uncountable family of continua no one of which can be continuously mapped onto any other. This was improved e.g. in \cite[20.3, 20.4, 20.9]{MackowiakTymchatyn} or \cite[page 2]{KrzempekPol} into a result stating that there is a family $\mathcal{Z}$ of size continuum of pairwise non-homeomorphic continua such that for every continuum $X$ there are at most countably many members of $\mathcal{Z}$ onto which $X$ can be continuously mapped.

\begin{lemma}\label{LemmaGraphCountableDegree}
Let $G$ be a directed graph with $\mathfrak c$-many vertices and with a countable outdegree at every vertex. Then there is a set of vertices $T$ of size $\mathfrak c$, such that for no pair of distinct vertices $u, v\in T$ both $(u,v)$ and $(v,u)$ form an edge in the graph.
\end{lemma}

\begin{proof}
Using transfinite recursion, let us construct a sequence $v_\alpha$, $\alpha<\mathfrak c$, of pairwise distinct vertices such that $(v_\alpha, v_\beta)$ is not an edge when $\alpha<\beta<\mathfrak c$. Given $\alpha<\mathfrak c$, assume we have already constructed $v_\gamma$ for each $\gamma<\alpha$. Let $M$ be the set of vertices $v$ in $G$ for which there is $\gamma<\alpha$ such that $(v_\gamma, v)$ is an edge. Clearly, $|M|\leq |\alpha|\cdot\omega<\mathfrak c$. Hence, there is a vertex $v_\alpha$ in $G$ which is not in $M \cup \{ v_{\gamma} \, ; \ \gamma < \alpha \}$. It follows that $(v_\gamma,v_\alpha)$ is not an edge in $G$ for any $\gamma<\alpha$.
\end{proof}

\begin{theorem}\label{ThmInfiniteDimCase}
There is a family $\mathcal{C}$ of size continuum of pairwise non-homeomorphic continua such that $X \times X$ is homeomorphic to $X$ for each $X \in \mathcal{C}$.
\end{theorem}

\begin{proof}
Consider the family $\mathcal{Z}$ described in the first paragraph of this section. Let $G$ be a directed graph with $\mathcal{Z}$ as the set of vertices such that $(X, Y) \in \mathcal{Z} \times \mathcal{Z}$ forms an edge in $G$ if and only if there is a continuous mapping from $X^{\mathbb N}$ onto $Y^{\mathbb N}$. Then $G$ has $\mathfrak c$-many vertices and it has a countable outdegree at every vertex. By Lemma \ref{LemmaGraphCountableDegree}, there is a family $\mathcal{Z}_0 \subseteq \mathcal{Z}$ of size $\mathfrak{c}$ such that for every two distinct continua $X,Y \in \mathcal{Z}_0$ either $X^{\mathbb N}$ can not be continuously mapped onto $Y^{\mathbb N}$ or $Y^{\mathbb N}$ can not be continuously mapped onto $X^{\mathbb N}$ (either way, $X^{\mathbb N}$ is not homeomorphic to $Y^{\mathbb N}$). Letting $\mathcal{C}:= \{ X^{\mathbb N} \, ; \ X \in \mathcal{Z}_0 \}$ completes the proof.
\end{proof}

\begin{remark}
An alternative proof of Theorem \ref{ThmInfiniteDimCase} can be obtained using Cook continua. Let us recall that Cook constructed a non-degenerate (even hereditarily indecomposable) continuum with no continuous mappings between its subcontinua except for constant and identity mappings \cite{Cook}. Any such continuum is usually called a Cook continuum. Later Mackowiak constructed a hereditary decomposable arc-like Cook continuum $C$ with no cut-points \cite[Theorem 6.1]{Mackowiak}. By the theory of tranches (see e.g. \cite[p. 200]{Kuratowski2}) there exists a monotone surjective mapping $f: C\to [0,1]$ (possessing a decomposition into tranches). Letting $C_s := f^{-1} (\{ s \})$ for every $s \in (0,1)$, the family $\{ C_s \, ; \, s \in (0,1) \}$ is of cardinality $\mathfrak c$ and it consists of pairwise disjoint nondegenerate (as $C$ has no cut-points) Cook continua. Moreover, for all $s,t \in (0,1)$, any continuous mapping between subcontinua of $C_s$ and $C_t$ is either constant, or it is the identity and $s=t$. Finally, given any two distinct numbers $s,t \in (0,1)$, assume towards contradiction that there is a homeomorphism $h \colon C_s^{\N}\to C_t^{\N}$. Then, letting $\Delta \colon C_s \to C_s^{\N}$ be the diagonal mapping, it follows that $\pi_i \circ h \circ \Delta: C_s \to C_t$ is a constant mapping for every $i \in \N$. Thus, $h \circ \Delta$ is a constant mapping, which is a contradiction.
\end{remark}

\begin{remark}
Yet another proof of Theorem \ref{ThmInfiniteDimCase} can be obtained from the following general result by Orsatti and Rodin\`o \cite{OrsattiRodino}. They have shown that for every $r \in \N$ and every infinite cardinal number $\lambda$, there is a class $\mathcal{C}$ of size $2^{\lambda}$ of pairwise non-homeomorphic compact connected Hausdorff topological Abelian groups of weight $\lambda$ with the property that for all $m,n \in \N$ and $X \in \mathcal{C}$, $X^m$ is homeomorphic to $X^n$ if and only if $m=n \ (\mathrm{mod} \ r)$. Hence, as every Hausdorff compact space of weight $\omega$ is metrizable, it suffices to take $\lambda=\omega$ and $r=1$.
\end{remark}

A question arises as to whether Theorem \ref{ThmInfiniteDimCase} remains true even if we additionally require each of the continua to be Peano. It turns out that such a strengthening of Theorem \ref{ThmInfiniteDimCase} can be realized using algebraic topology invariants.

\begin{theorem}
There is a family $\mathcal{P}$ of size continuum of pairwise non-homeomorphic Peano continua such that $X \times X$ is homeomorphic to $X$ for each $X \in \mathcal{P}$.
\end{theorem}

\begin{proof}
Let $\mathbb P$ be the set of primes.
For every $p\in\mathbb P$ fix a Peano continuum $Y_p$ whose fundamental group is isomorphic to $\mathbb Z_p$ (it is a folklore result that every finitely presented group can be realized as the fundamental group of a compact, connected, smooth manifold of dimension 4). For every nonempty set $A \subseteq \mathbb{P}$, let $X_A:=\prod_{p\in A} Y_p$. The fundamental group of $X_A$ is isomorphic to $G_A:=\prod_{p\in A} \mathbb{Z}_p$ since 
by \cite[Proposition 4.2]{Hatcher} the fundamental group of a product (even infinite) is isomorphic to the corresponding product of fundamental groups.
Consequently $X_A$ is not homeomorphic to $X_B$ if $A \neq B$, since the groups $G_A$ and $G_B$ are not isomorphic (if $p\in A\setminus B$, then $G_A$ contains a point of order $p$, whereas $G_B$ does not contain a point of order $p$).
\end{proof}

In the preceding proof we used the tool of infinite powers to obtain spaces homeomorphic to their respective squares and the algebraic tool of fundamental groups to prove that the spaces are pairwise non-homeomorphic. However, there are natural classes of (Peano) continua where neither of these tools can be used, e.g. countable-dimensional continua (see \cite[5.1]{EngelkingDimension}) or continua with trivial shape (see \cite{DydakSegal}).

\begin{question}
Is there a non-degenerate countable-dimensional (Peano) continuum $X$ homeomorphic to its square? If so, how many such continua are there (up to homeomorphism)?
\end{question}

Note that the non-existence of a non-degenerate countable-dimensional continuum homeomorphic to its square would immediately follow if we knew that $\mathrm{trind}(X)<\mathrm{trind}(X^2)$ for every infinite-dimensional countable-dimensional metrizable compact space $X$ (see \cite[Corollary 7.1.32]{EngelkingDimension}).

\begin{question}
How many (Peano) continua with trivial shape homeomorphic to their respective squares are there (up to homeomorphism)?
\end{question}

Going further, an absolute retract homeomorphic to its square is either trivial or it is homeomorphic to the Hilbert cube \cite{vanMillExample}. Note that, as was shown by Borsuk \cite[Corollary 11.2]{Borsuk}, there is a locally contractible and contractible (hence Peano) continuum which is not an absolute retract.

\begin{question}
How many contractible (Peano) continua homeomorphic to their respective squares are there (up to homeomorphism)?
\end{question}

\section{Zero-dimensional case}\label{SectionZeroDim}
Countable compact spaces will play an important role in this section. The techniques used to work with such spaces date back to Mazurkiewicz and Sierpi\'nski. In the first volume of Fundamenta Mathematicae \cite{MazurkiewiczSierpinski} they showed that every infinite countable compact metrizable space $X$ is homeomorphic to a countable ordinal of the form $\omega^\alpha \cdot k+1$, where $\alpha = \mathrm{CB}(X)>0$ and $k \in \omega \setminus \{ 0 \}$ is the cardinality of $X^{(\alpha)}$.

For each $\alpha < \omega_1$, consider the space $Z(\alpha)$ from \cite{CharatonikSahan}. That is, $Z(\alpha)$ is uncountable, metrizable, compact, zero-dimensional and it satisfies the following two conditions:
\begin{itemize}
    \item $\mathrm{CB} (x,Z(\alpha)) < \alpha$ for every $x \in Z(\alpha) \setminus \mathrm{PK} (Z(\alpha))$;
    \item $\mathrm{PK} (Z(\alpha)) \subseteq \overline{\big\lbrace x \in Z(\alpha) \setminus \mathrm{PK} (Z(\alpha)) \, ; \ \mathrm{CB} (x,Z(\alpha)) = \beta \big\rbrace}$ for every $\beta < \alpha$.
\end{itemize}

The following two lemmata follow from Proposition \ref{IsolatedPointsAreDenseInTheComplementOfThePerfectKernel} and \cite[Theorem 3.4]{CharatonikSahan}, respectively.

\begin{lemma}\label{IsolatedPointsAreDenseInZ(alpha)}
    For every $\alpha < \omega_1$ with $\alpha \neq 0$, the set $Z(\alpha) \setminus (Z(\alpha))'$ is dense in $Z(\alpha)$.
\end{lemma}

\begin{lemma}\label{UnctblClopenIsHomeo}
    Let $\alpha < \omega_1$. Then every uncountable clopen subset of $Z(\alpha)$ is homeomorphic to $Z(\alpha)$.
\end{lemma}

Let $\mathcal{O} := \{ 1 \} \cup \{ \omega^n +1 \, ; \ n \in \omega \, , \ n \neq 0 \}$. Clearly, every ordinal in $\mathcal{O}$ is (with the order topology) a countable metrizable compact space with finite Cantor-Bendixson rank.

Now we are ready to describe the basic strategy behind the proof of Theorem \ref{MainResult}. For each infinite subset $M$ of $\omega \setminus \{ 0 \}$, we will consider the family $\mathscr{S}(M)$ of all finite products of members of $\mathcal{O} \cup \{ Z(m) \, ; \, m \in M \}$. We will prove that every $S \in \mathscr{S}(M)$ admits a finite $S$-clopen partition into arbitrarily small sets such that each member of the partition is homeomorphic to a member of $\mathscr{S}(M)$. This will allow us to prove the uniqueness of a specific compactification $X(M)$ of the topological sum of the family $\mathscr{S}(M)$, where each $S \in \mathscr{S}(M)$ is taken infinitely (countably) many times in the sum. The space $X(M)$ will be constructed in such a way that its topological characterization will make it possible to verify (with the help of Lemma \ref{UnctblClopenIsHomeo}) that $X(M) \times X(M)$ is homeomorphic to $X(M)$ and that $X(M)$ is not homeomorphic to $X(L)$ for any other infinite set $L \subseteq \omega \setminus \{ 0 \}$.

\begin{lemma}\label{ClopenPartition1}
    Let $X$ be a nonempty countable compact metric space with finite Cantor-Bendixson rank and let $\varepsilon>0$. Then there is a finite $X$-clopen partition $\mathcal{F}$ of $X$ with $\mathrm{mesh}(\mathcal{F})<\varepsilon$ such that every $F \in \mathcal{F}$ is homeomorphic to a member of $\mathcal{O}$. If, moreover, $X$ is homeomorphic to a member of $\mathcal{O}$, then at least one element of $\mathcal{F}$ is homeomorphic to $X$.
\end{lemma}
\begin{proof}
    Since $X$ is compact and zero-dimensional, there is a finite $X$-clopen partition $\mathcal{A}$ of $X$ with $\mathrm{mesh}(\mathcal{A})<\varepsilon$. We claim that, for every $A \in \mathcal{A}$, there is a finite $X$-clopen partition $\mathcal{F}_A$ of $A$ such that every element of $\mathcal{F}_A$ is homeomorphic to a member of $\mathcal{O}$. Let $A \in \mathcal{A}$ be given. If $A$ is finite, we define $\mathcal{F}_A := \big\lbrace \{ x \} ; \, x \in A \big\rbrace$. If $A$ is infinite, then, as $\mathrm{CB(A)} \leq \mathrm{CB}(X) < \omega$, there are $k,n \in \omega \setminus \{ 0 \}$ such that $A$ is homeomorphic to $\omega^n \cdot k+1$. In that case, however, there is an $A$-clopen (and thus also $X$-clopen) partition $\mathcal{F}_A$ of $A$ with exactly $k$ elements each of which is homeomorphic to $\omega^n+1$. Let $\mathcal{F}:= \bigcup \{ \mathcal{F}_A \, ; \, A \in \mathcal{A} \}$.

    Assume that $X$ is homeomorphic to $\omega^n +1$ for some $n \in \omega \setminus \{0\}$. Then $X^{(n)}=\{ x \}$ for some $x \in X$. Consequently, for any $X$-clopen partition $\mathcal{F}$ of $X$, the member $F$ of $\mathcal{F}$ containing $x$ satisfies $F^{(n)}=\{ x \}$ and thus is homeomorphic to $\omega^n +1$.
\end{proof}

\begin{lemma}\label{ClopenPartition2}
    Let $n \in \omega$ and let $X$ be a metric space homeomorphic to $Z(n)$. For every $\varepsilon >0$, there is a finite $X$-clopen partition $\mathcal{F}$ of $X$ with $\mathrm{mesh}(\mathcal{F})<\varepsilon$ such that every $F \in \mathcal{F}$ is homeomorphic either to $Z(n)$ or to a member of $\mathcal{O}$. In particular, at least one member of $\mathcal{F}$ is homeomorphic to $Z(n)$.
\end{lemma}
\begin{proof}
    Since $X$ is compact and zero-dimensional, there is a finite $X$-clopen partition $\mathcal{A}$ of $X$ with $\mathrm{mesh}(\mathcal{A})<\varepsilon$. Let $\mathcal{A}_0:=\{ A \in \mathcal{A} \, ; \ A \textup{ is countable} \}$ and $\mathcal{A}_1 := \mathcal{A} \setminus \mathcal{A}_0$. By Lemma \ref{UnctblClopenIsHomeo}, every member of $\mathcal{A}_1$ is homeomorphic to $Z(n)$. For every $A \in \mathcal{A}_0$, since $n<\omega$, we have $\mathrm{CB}(A)<\omega$. Hence, by Lemma \ref{ClopenPartition1}, there is a finite $A$-clopen (and thus also $X$-clopen) partition $\mathcal{F}_A$ of $A$ with $\mathrm{mesh}(\mathcal{F}_A)<\varepsilon$ such that every $F \in \mathcal{F}_A$ is homeomorphic to a member of $\mathcal{O}$. Let $\mathcal{F}:=\mathcal{A}_1 \cup \bigcup \{ \mathcal{F}_A \, ; \, A \in \mathcal{A}_0 \}$. 
\end{proof}

For every infinite set $M \subseteq \omega$, let $\mathscr{S}(M)$ be a family of topological spaces defined by
\[ \mathscr{S}(M) := \big\lbrace \alpha_0 \times \dots \times \alpha_k \times Z(m_1) \times \dots \times Z(m_n) \, ; \ k,n \in \omega, \ \alpha_0 , \dotsc ,\alpha_k \in \mathcal{O}, \ m_1, \dotsc ,m_n \in M \big\rbrace . \]
Clearly, $\mathscr{S}(M)$ is countable and every member of $\mathscr{S}(M)$ is a nonempty compact metrizable zero-dimensional space. Using Lemmata \ref{ClopenPartition1} and \ref{ClopenPartition2} we easily deduce the following.

\begin{lemma}\label{ClopenPartition3}
    Assume $M \subseteq \omega$ is an infinite set, let $Y$ be a metric space homeomorphic to a member of $\mathscr{S}(M)$ and let $\varepsilon>0$. Then there is a finite $Y$-clopen partition $\mathcal{F}$ of $Y$ with $\mathrm{mesh}(\mathcal{F})<\varepsilon$ such that every $F \in \mathcal{F}$ is homeomorphic to a member of $\mathscr{S}(M)$ and at least one element of $\mathcal{F}$ is homeomorphic to $Y$.
\end{lemma}

As a consequence we get the following lemma.

\begin{lemma}\label{ClopenPartition4}
    Let $M \subseteq \omega$ be an infinite set, let $X$ be a compact metric space and let $\mathcal{A}$ be a family of pairwise disjoint clopen subsets of $X$ such that every member of $\mathcal{A}$ is homeomorphic to a member of $\mathscr{S}(M)$. Then there is a null family $\mathcal{A}^*$ of pairwise disjoint clopen subsets of $X$ such that:
    \begin{enumerate}[label=(\roman*),font=\textup,noitemsep]
        \item $\bigcup \mathcal{A}^* = \bigcup \mathcal{A}$;
        \item for every $B \in \mathcal{A}^*$, there is $A \in \mathcal{A}$ with $B \subseteq A$;
        \item for every $A \in \mathcal{A}$, there is $B \in \mathcal{A}^*$ with $B \subseteq A$ such that $A$ is homeomorphic to $B$;
        \item every member of $\mathcal{A}^*$ is homeomorphic to a member of $\mathscr{S}(M)$.
    \end{enumerate}
\end{lemma}
\begin{proof}
    If $\mathcal{A}$ is finite, we just let $\mathcal{A}^* := \mathcal{A}$. Assume $\mathcal{A}$ is infinite. Since $X$ is separable, $\mathcal{A}$ is countable. Hence, we can write $\mathcal{A}= \{ Y_n \, ; \ n \in \omega \}$, where $Y_k \neq Y_n$ (and thus $Y_k \cap Y_n = \emptyset$) when $k \neq n$. For every $n \in \omega$, there is (by Lemma \ref{ClopenPartition3}) a finite $Y_n$-clopen partition $\mathcal{F}_n$ of $Y_n$ with $\mathrm{mesh}(\mathcal{F}_n)< 2^{-n}$ such that every $F \in \mathcal{F}_n$ is homeomorphic to a member of $\mathscr{S}(M)$ and at least one element of $\mathcal{F}_n$ is homeomorphic to $Y_n$. Let $\mathcal{A}^* := \bigcup \{ \mathcal{F}_n \, ; \ n \in \omega \}$.
\end{proof}

Several variants of the next proposition are well known and were reproved many times for different purposes (e.g. \cite{KnasterReichbach, Pelczynski, Lorch, Terasawa, Zielinski},  or \cite[Proposition 8.8]{IllanesNadler}).
For the sake of completeness we include the proof.

\begin{proposition}\label{Colors}
    Let $X_1$ and $X_2$ be infinite compact metrizable spaces and let $M$ be a nonempty countable set. For each $i \in \{ 1,2 \}$, let $\mu_i \colon X_i \setminus X_i' \to M$ be a mapping such that, for every $m \in M$,
    \[ X_i' \subseteq \overline{ \big\lbrace x \in X_i \setminus X_i' \, ; \ \mu_i (x) = m \big\rbrace }. \tag{$\ast$} \]
    Let $h \colon X_1' \to X_2'$ be a homeomorphism. Then there is a homeomorphism $\overline{h} \colon X_1 \to X_2$ extending $h$ such that $\mu_1 (x) = \mu_2 \big( \overline{h}(x) \big)$ for every $x \in X_1 \setminus X_1'$.
\end{proposition}

\begin{proof}
    Since $X_1$, $X_2$ are infinite and compact, the sets $X_1'$, $X_2'$ are nonempty. Hence, it follows from $(\ast)$ that $X_1 \setminus X_1'$ and $X_2 \setminus X_2'$ are infinite sets. On the other hand, these two sets are countable as $X_1$, $X_2$ are separable. Thus, for each $i \in \{ 1,2 \}$, we can fix a bijection $\varphi_i \colon \omega \to X_i \setminus X_i'$. Let $\varrho$ and $\sigma$ be compatible metrics on $X_1$ and $X_2$, respectively. Using $(\ast)$, it is easy to construct (inductively) sequences $(a_n)_{n \in \omega}$, $(b_n)_{n \in \omega}$ of elements of $\omega$ in such a way that the following holds for each $n \in \omega$:
    \begin{itemize}
        \item $\mu_1 \big( \varphi_1 (a_n) \big) = \mu_2 \big( \varphi_2 (b_n) \big)$.
        \item If $n$ is even, then $b_n \notin \{ b_i \, ; \ i<n \}$, $a_n = \min \big( \omega \setminus \{ a_i \, ; \ i<n \} \big)$ and there is $p \in X_1'$ such that $\varrho \big(\varphi_1(a_n),p \big) = \mathrm{dist}_{\varrho} \big( \varphi_1(a_n), X_1' \big)$ and $\sigma \big( \varphi_2 (b_n), h(p) \big) < 2^{-n}$.
        \item If $n$ is odd, then $a_n \notin \{ a_i \, ; \ i<n \}$, $b_n = \min \big( \omega \setminus \{ b_i \, ; \ i<n \} \big)$ and there is $q \in X_2'$ such that $\sigma \big(\varphi_2(b_n),q \big) = \mathrm{dist}_{\sigma} \big( \varphi_2 (b_n), X_2' \big)$ and $\varrho \big( \varphi_1 (a_n), h^{-1}(q) \big) < 2^{-n}$.
    \end{itemize}
    Clearly, $\{ a_n \, ; \ n \in \omega \} = \{ b_n \, ; \ n \in \omega \} = \omega$ and $a_i \neq a_j$, $b_i \neq b_j$ for all $i,j \in \omega$ with $i \neq j$. Therefore, the mapping $g \colon X_1 \setminus X_1' \to X_2 \setminus X_2'$ given by $g \big( \varphi_1 (a_n) \big) = \varphi_2 (b_n)$, $n \in \omega$, is a well-defined bijection. Define a mapping $\overline{h} \colon X_1 \to X_2$ by $\overline{h}(x)=h(x)$ for $x \in X_1'$ and by $\overline{h}(x)=g(x)$ for $x \in X_1 \setminus X_1'$. It is clear that $\mu_1 (x) = \mu_2 \big( \overline{h}(x) \big)$ for every $x \in X_1 \setminus X_1'$, it remains to show that $\overline{h}$ is a homeomorphism. Since $\overline{h}$ is a bijection and $X_1$ is compact, it suffices to prove that $\overline{h}$ is continuous. Trivially, $\overline{h}$ is continuous at every point of the set $X_1 \setminus X_1'$. Given $z \in X_1'$ and $\varepsilon > 0$, let us find $\delta>0$ such that $\sigma \big( \overline{h}(z), \overline{h}(x) \big) < 2 \varepsilon$ for every $x \in X_1$ with $\varrho (z,x)<\delta$. By the continuity of $h$, there is $\delta_1 >0$ such that $\sigma \big( h(z), h(x) \big) < \varepsilon$ for every $x \in X_1'$ with $\varrho (z,x)< 2\delta_1$. Also, there is $\delta_2 >0$ such that $2^{-n} < \min \{ \varepsilon, \delta_1 \}$ for every $n \in \omega$ with $\varrho \big( z, \varphi_1 (a_n) \big) < \delta_2$. By the compactness of $X_2$, there is $\delta_3 >0$ such that $\mathrm{dist}_{\sigma} \big( \varphi_2 (b_n), X_2' \big) < \varepsilon$ for every $n \in \omega$ with $\varrho \big( z, \varphi_1 (a_n) \big) < \delta_3$. Let $\delta:= \min \{ \delta_1, \delta_2, \delta_3 \}$ and let $x \in X_1$ satisfy $\varrho (z,x)<\delta$. If $x \in X_1'$, then (as $\delta \leq \delta_1$) we immediately receive $\sigma \big( \overline{h}(z), \overline{h}(x) \big) = \sigma \big( h(z), h(x) \big) < \varepsilon$. Assume $x \in X_1 \setminus X_1'$ and let $n \in \omega$ satisfy $\varphi_1 (a_n) = x$. Since $\delta \leq \min \{ \delta_2, \delta_3 \}$, we have $2^{-n} < \min \{ \varepsilon, \delta_1 \}$ and $\mathrm{dist}_{\sigma} \big( \varphi_2 (b_n), X_2' \big) < \varepsilon$. If $n$ is even, there is $p \in X_1'$ such that $\sigma \big( \varphi_2 (b_n), h(p) \big) < 2^{-n} < \varepsilon$ and
    \[ \varrho (x,p) = \varrho \big(\varphi_1(a_n),p \big) = \mathrm{dist}_{\varrho} \big( \varphi_1(a_n), X_1' \big) \leq \varrho \big( \varphi_1(a_n), z \big) = \varrho (z,x). \]
    Then $\varrho (z,p) \leq \varrho (z,x) + \varrho (x,p) \leq 2\varrho (z,x) < 2\delta \leq 2\delta_1$. Hence, $\sigma \big( h(z), h(p) \big) < \varepsilon$, implying that
    \[  \sigma \big( \overline{h}(z), \overline{h}(x) \big) \leq \sigma \big( \overline{h}(z), \overline{h}(p) \big) + \sigma \big( \overline{h}(p), \overline{h}(x) \big) = \sigma \big( h(z), h(p) \big) + \sigma \big( h(p), \varphi_2(b_n) \big) < 2 \varepsilon . \]
    If $n$ is odd, there is $q \in X_2'$ such that $\varrho \big( x, h^{-1}(q) \big) = \varrho \big( \varphi_1 (a_n), h^{-1}(q) \big) < 2^{-n} < \delta_1$ and
    \[ \sigma \big(\overline{h}(x),q \big) = \sigma \big( \varphi_2(b_n),q \big) = \mathrm{dist}_{\sigma} \big( \varphi_2(b_n), X_2' \big) < \varepsilon. \]
    Then, denoting $p:=h^{-1}(q)$, we have $\varrho (z,p) \leq \varrho (z,x) + \varrho (x,p) < \delta + \delta_1 \leq 2\delta_1$. Therefore, we obtain $\sigma \big( \overline{h}(z), q \big) = \sigma \big( h(z), h(p) \big) < \varepsilon$. Thus, $\sigma \big( \overline{h}(z), \overline{h}(x) \big) \leq \sigma \big( \overline{h}(z), q \big) + \sigma \big( q, \overline{h}(x) \big) < 2 \varepsilon$.
\end{proof}

\begin{proposition}\label{ExtendingHomeos}
    Let $\mathscr{S}$ be a countable family of nonempty metrizable compact topological spaces and let $(X_1, d^1)$, $(X_2, d^2)$ be compact metric spaces. For each $i \in \{ 1,2 \}$, let $\mathcal{A}_i$ be a null family of pairwise disjoint clopen subsets of $X_i$ such that:
    \begin{enumerate}[label=(\arabic*),font=\textup]
        \item $C_i := X_i \setminus \bigcup \mathcal{A}_i$ is a perfect set;
        \item every member of $\mathcal{A}_i$ is homeomorphic to a member of $\mathscr{S}$;
        \item for every $S \in \mathscr{S}$ and every open set $V \subseteq X_i$ with $V \cap C_i \neq \emptyset$, there is $A \in \mathcal{A}_i$ such that $A$ is homeomorphic to $S$ and $A \subseteq V$.
    \end{enumerate}
    Then every homeomorphism $h \colon C_1 \to C_2$ can be extended to a homeomorphism $\overline{h} \colon X_1 \to X_2$.
\end{proposition}
\begin{proof}
    Let $h \colon C_1 \to C_2$ be a homeomorphism. For each $i \in \{ 1,2 \}$, define $\mathcal{D}_i := \mathcal{A}_i \cup \big\lbrace \{x\} ; \, x \in C_i \big\rbrace$. Clearly, $\mathcal{D}_i \subseteq \mathcal{K}(X_i)$.
    
    \begin{claim}\label{ClosednessOfPartition2}
        For each $i \in \{ 1,2 \}$, the set $\mathcal{D}_i$ is closed in $\mathcal{K}(X_i)$.
    \end{claim}
    \begin{claimproof}
        Given $i \in \{ 1,2 \}$ and $K \in \mathcal{K}(X_i) \setminus \mathcal{D}_i$, let us find an open set $\mathcal{U} \subseteq \mathcal{K}(X_i)$ such that $K \in \mathcal{U}$ and $\mathcal{U} \cap \mathcal{D}_i = \emptyset$. If $K \cap A \neq \emptyset$ for some $A \in \mathcal{A}_i$, we just let $\mathcal{U}:= \big\lbrace L \in \mathcal{K}(X_i) \, ; \ L \cap A \neq \emptyset \big\rbrace \setminus \{ A \}$. Thus, assume $K \subseteq C_i$. As $K \notin \mathcal{D}_i$, the set $K$ is not a singleton. Hence, $r:= \mathrm{diam}_{d^i} (K)$ is positive. Since $\mathcal{A}_i$ is a null family, $\mathcal{F} := \big\lbrace A \in \mathcal{A}_i \, ; \ \mathrm{diam}_{d^i} (A) > r/2 \big\rbrace$ is finite and thus closed in $\mathcal{K}(X_i)$. Let $\mathcal{U}:=\big\lbrace L \in \mathcal{K}(X_i) \, ; \ \mathrm{diam}_{d^i}(L)>r/2 \big\rbrace \setminus \mathcal{F}$. Then $K \in \mathcal{U} \subseteq \mathcal{K}(X_i) \setminus \mathcal{D}_i$ and $\mathcal{U}$ is open in $\mathcal{K}(X_i)$. \claimend
    \end{claimproof}
    
    As $X_1$, $X_2$ are compact, so are $\mathcal{K}(X_1)$, $\mathcal{K}(X_2)$. Hence, $\mathcal{D}_1$ and $\mathcal{D}_2$ are compact by Claim \ref{ClosednessOfPartition2}. For each $i \in \{ 1,2 \}$, since every member of $\mathcal{A}_i$ is an isolated point of $\mathcal{D}_i$ and since $\big\lbrace \{x\} ; \, x \in C_i \big\rbrace$ is homeomorphic to the perfect set $C_i$, we have $\mathcal{D}_i' = \big\lbrace \{x\} ; \, x \in C_i \big\rbrace$. Let $\mathscr{S}_0$ be a subfamily of $\mathscr{S}$ such that for every $S \in \mathscr{S}$ there is exactly one $T \in \mathscr{S}_0$ homeomorphic to $S$. Then we have by (2) that, for each $i \in \{ 1,2 \}$, there is a unique mapping $\mu_i \colon \mathcal{A}_i \to \mathscr{S}_0$ such that $\mu_i (A)$ is homeomorphic to $A$ for every $A \in \mathcal{A}_i$. For each $i \in \{ 1,2 \}$ and $S \in \mathscr{S}_0$, it follows from (3) that
    \[ \mathcal{D}_i' \subseteq \overline{\big\lbrace A \in \mathcal{A}_i \, ; \ A \textup{ is homeomorphic to } S \big\rbrace} = \overline{\big\lbrace A \in \mathcal{D}_i \setminus \mathcal{D}_i' \, ; \ \mu_i(A)=S \big\rbrace} \]
    in $\mathcal{K}(X_i)$. Define a mapping $\varphi \colon \mathcal{D}_1' \to \mathcal{D}_2'$ by $\varphi (\{ x \}) = \{ h(x) \}$, $x \in C_1$. Since $h$ is a homeomorphism, so is $\varphi$. By Proposition \ref{Colors}, there is a homeomorphism $\overline{\varphi} \colon \mathcal{D}_1 \to \mathcal{D}_2$ which extends $\varphi$ and satisfies $\mu_1(A)=\mu_2 \big( \overline{\varphi}(A) \big)$ for every $A \in \mathcal{A}_1$. Fix a homeomorphism $h_A \colon A \to \overline{\varphi}(A)$ for every $A \in \mathcal{A}_1$ and define a mapping $\overline{h} \colon X_1 \to X_2$ by $h(x)=x$ for $x \in C_1$ and by $h(x)=h_A(x)$ for $x \in A \in \mathcal{A}_1$. Then $\overline{h}$ is a well-defined bijection and, since $\mathcal{A}_1$ is an $X_1$-clopen partition of $X_1 \setminus C_1$, it is continuous at every point of $X_1 \setminus C_1$. Given $x \in C_1$ and $\varepsilon > 0$, let us find $\delta > 0$ such that $d^2 \big( \overline{h}(x), \overline{h}(y) \big) < \varepsilon$ for every $y \in X_1$ with $d^1(x,y)<\delta$. By the continuity of $\overline{\varphi}$, there is $\delta_1 > 0$ such that $d_H^2 \big( \overline{\varphi}(D), \overline{\varphi} (\{ x \}) \big) < \varepsilon$ for every $D \in \mathcal{D}_1$ with $d_H^1 (D,\{ x \}) < 2 \delta_1$. Since $\mathcal{A}_1$ is a null family, so is $\mathcal{D}_1$. Hence, as $\mathcal{D}_1$ consists of closed sets and $\{ x \}$ is the only member of $\mathcal{D}_1$ containing $x$, there is $\delta_2 > 0$ such that $\mathrm{diam}_{d^1} (D) < \delta_1$ for every $D \in \mathcal{D}_1$ with $\mathrm{dist}_{d^1}(x,D)< \delta_2$. Let $\delta := \min \{ \delta_1, \delta_2 \}$ and let $y \in X_1$ be any point satisfying $d^1 (x,y)< \delta$. There is $D \in \mathcal{D}_1$ with $y \in D$. Clearly, $\mathrm{dist}_{d^1}(x,D) \leq d^1(x,y)< \delta
    \leq \delta_2$, hence $\mathrm{diam}_{d^1} (D) < \delta_1$. Thus, $d_H^1 (D,\{ x \}) = \max \big\lbrace d^1 (x,z) \, ; \ z \in D \big\rbrace \leq d^1(x,y) + \mathrm{diam}_{d^1} (D) < \delta + \delta_1 \leq 2 \delta_1$, which gives us $d_H^2 \big( \overline{\varphi}(D), \overline{\varphi} (\{ x \}) \big) < \varepsilon$. However, by the definition of $\overline{h}$, we have $\overline{h} (y) \in \overline{\varphi}(D)$ and $\big\lbrace \overline{h}(x) \big\rbrace = \overline{\varphi} (\{ x \})$, hence $d^2 \big( \overline{h}(x), \overline{h}(y) \big) \leq d_H^2 \big( \big\lbrace \overline{h}(x) \big\rbrace , \overline{\varphi}(D) \big) = d_H^2 \big( \overline{\varphi}(D), \overline{\varphi} (\{ x \}) \big) < \varepsilon$.

    Having shown that $\overline{h}$ is continuous, it follows from the compactness of $X_1$ that $\overline{h}$ is a homeomorphism.
\end{proof}

\begin{proposition}\label{ExistenceAndUniqenessOfTheCompactification}
    Let $M \subseteq \omega$ be an infinite set. Then there exists a compact metrizable space $X$ with the property that there is a family $\mathcal{A}$ of pairwise disjoint clopen subsets of $X$ such that:
    \begin{enumerate}[label=(\arabic*),font=\textup]
        \item $C := X \setminus \bigcup \mathcal{A}$ is homeomorphic to the Cantor space;
        \item every member of $\mathcal{A}$ is homeomorphic to a member of $\mathscr{S}(M)$;
        \item for every $S \in \mathscr{S}(M)$ and every open set $V \subseteq X$ with $V \cap C \neq \emptyset$, there is $A \in \mathcal{A}$ such that $A$ is homeomorphic to $S$ and $A \subseteq V$.
    \end{enumerate}
    Moreover, the space $X$ is unique up to homeomorphism.
\end{proposition}

\begin{proof}
    Let $(S_n)_{n \in \omega}$ be a sequence of spaces in $\mathscr{S}(M)$ such that $\{ n \in \omega \, ; \ S_n = S \}$ is an infinite set for every $S \in \mathscr{S}(M)$. Fix a set $C_0 \subseteq \R$ homeomorphic to the Cantor space and let $\{ q_k \, ; \ k \in \omega \}$ be a countable dense subset $C_0$. Let $(t_n)_{n \in \omega}$ be a strictly decreasing sequence of positive real numbers converging to $0$. For every $n \in \omega$, since $S_n$ is separable, metrizable and zero-dimensional, it can be embedded in $\R$. Thus, for each $n \in \omega$, there is a set $A_n \subseteq (t_{n+1},t_{n})$ homeomorphic to $S_n$. Define $\mathcal{A}:= \big\lbrace \{ q_k \} \times A_n \, ; \ k,n \in \omega \, , \ k \leq n \big\rbrace$ and $X := (C_0 \times \{ 0 \} ) \cup \bigcup \mathcal{A}$. It is easy to see that $X$ is a compact subset of $\R^2$ and that $\mathcal{A}$ is a disjoint family consisting of relatively clopen subsets of $X$. Obviously, (2) is satisfied and the set $C := X \setminus \bigcup \mathcal{A} = C_0 \times \{ 0 \}$ is homeomorphic to the Cantor space. Now, given $S \in \mathscr{S}(M)$ and an open set $V \subseteq \R^2$ with $V \cap C \neq \emptyset$, let us show that there are $k,n \in \omega$ with $k \leq n$ such that $\{ q_k \} \times A_n$ is a subset of $V$ homeomorphic to $S$. Since $V \cap C \neq \emptyset$, there is $k \in \omega$ with $(q_k,0) \in V$. Since $V$ is open, there is $\varepsilon > 0$ such that $\{ q_k \} \times (0, \varepsilon) \subseteq V$. Fix $m \in \omega$ with $t_m \leq \varepsilon$. Since $\{ n \in \omega \, ; \ S_n = S \}$ is infinite, there is $n \in \omega$ such that $n \geq \max \{ k,m \}$ and $S_n = S$. Then $\{ q_k \} \times A_n$ is homeomorphic to $S$ and, as $A_n \subseteq (t_{n+1},t_n) \subseteq (0,t_m) \subseteq (0,\varepsilon)$, it is a subset of $V$.

    The uniqueness of $X$ easily follows from Lemma \ref{ClopenPartition4} and Proposition \ref{ExtendingHomeos}.
\end{proof}

For every infinite set $M \subseteq \omega$, denote by $X(M)$ the corresponding space whose existence and uniqueness was proven in Proposition \ref{ExistenceAndUniqenessOfTheCompactification}. Note that $X(M)$ is zero-dimensional since it is the union of countably many closed zero-dimensional sets.

\begin{proposition}\label{X(M)IsHomeoToItsSquare}
    Suppose $M \subseteq \omega \setminus \{ 0 \}$ is an infinite set. Then $X(M) \times X(M)$ is homeomorphic to $X(M)$.
\end{proposition}
\begin{proof}
    Let us write $X$ instead of $X(M)$ throughout this proof. Fixing a family $\mathcal{A}$ witnessing the defining property of the space $X$ (see the statement of Proposition \ref{ExistenceAndUniqenessOfTheCompactification}), denote $C := X \setminus \bigcup \mathcal{A}$. Since $\mathcal{A}$ is infinite and countable, we can write $\mathcal{A}= \{ A_i \, ; \, i \in \omega \}$, where $A_i \neq A_j$ (and thus $A_i \cap A_j = \emptyset$) for $i \neq j$. Let $\Gamma := \{ i \in \omega \, ; \, A_i \textup{ is not a singleton} \}$ and note that $A_i$ is infinite for every $i \in \Gamma$. Since $0 \notin M$, it follows from Lemma \ref{IsolatedPointsAreDenseInZ(alpha)} and Proposition \ref{IsolatedPointsAreDenseInTheComplementOfThePerfectKernel} that $A_i \setminus A_i'$ is dense in $A_i$ for every $i \in \omega$. In particular, $A_i \setminus A_i'$ is infinite for every $i \in \Gamma$. Hence, for each $i \in \Gamma$, there is an injective sequence $(a_{i,k})_{k \in \omega}$ with $A_i \setminus A_i' = \{ a_{i,k} \, ; \, k \in \omega \}$. For all $i \in \Gamma$ and $n \in \omega$, let
    \[ \mathcal{D}_i(n):= \big\lbrace \{ a_{i,k} \} \, ; \, 0 \leq k \leq n \big\rbrace \cup \big\lbrace A_i \setminus \{ a_{i,k} \, ; \, 0 \leq k \leq n \} \big\rbrace . \]
    In addition, let $\mathcal{D}_i(n):= \{ A_i \}$ for all $i \in \omega \setminus \Gamma$ and $n \in \omega$. Then, for all $i,n \in \omega$, the family $\mathcal{D}_i(n)$ is a finite $X$-clopen partition of $A_i$ and every member of $\mathcal{D}_i(n)$ is either a singleton or (by Lemma \ref{SubstractingOneIsolatedPointDoesntChangeAnything}) homeomorphic to $A_i$. In particular, every member of $\mathcal{D}_i(n)$ is homeomorphic to a member of $\mathscr{S}(M)$. Let $C^* := (C \times X) \cup (X \times C)$ and
    \[ \mathcal{A}^* :=  \bigcup_{i,j \in \omega} \big\lbrace E \times F \, ; \ E \in \mathcal{D}_i(i+j) , \, F \in \mathcal{D}_j(i+j) \big\rbrace . \]
    It easily follows from the classical topological characterization of the Cantor space due to Brouwer that $C^*$ is homeomorphic to the Cantor space. Moreover, we have
    \[ \bigcup \mathcal{A}^* = \bigcup_{i,j \in \omega} (A_i \times A_j) = \Big( \bigcup \mathcal{A} \Big) \times \Big( \bigcup \mathcal{A} \Big) = (X \setminus C) \times (X \setminus C) = (X \times X) \setminus C^* \]
    and it is clear that $\mathcal{A}^*$ is a family of pairwise disjoint clopen subsets of $X \times X$. Since it is obvious that the product of finitely many members of $\mathscr{S}(M)$ is homeomorphic to a member $\mathscr{S}(M)$, every member of $\mathcal{A}^*$ is homeomorphic to a member $\mathscr{S}(M)$. Given $S \in \mathscr{S}(M)$ and an open set $W \subseteq X \times X$ with $W \cap C^* \neq \emptyset$, let us show that there is $K \in \mathcal{A}^*$ with $K \subseteq W$ such that $K$ is homeomorphic to $S$. Let $U,V \subseteq X$ be open sets satisfying $U \times V \subseteq W$ and $(U \times V) \cap C^* \neq \emptyset$. Then $U \times V$ intersects at least one of the sets $C \times (X \setminus C)$, $(X \setminus C) \times C$ and $C \times C$. We will assume the second possibility is true (the first two possibilities are symmetric and the third one is easier to deal with), that is, $U \setminus C \neq \emptyset$ and $V \cap C \neq \emptyset$. Then there is $i \in \omega$ such that $U \cap A_i \neq \emptyset$. Let us assume that $i \in \Gamma$ (the situation is easier if $i \notin \Gamma$). Then, since $\{ a_{i,k} \, ; \, k \in \omega \}$ is dense in $A_i$, there is $k \in \omega$ with $a_{i,k} \in U$. Let $V_0 := V \setminus \bigcup \{ A_n \, ; \ 0 \leq n \leq k \}$, then $V_0$ is an open subset of $X$ intersecting $C$. Hence, as $\mathcal{A}$ witnesses the defining property of the space $X$, there is $j \in \omega$ such that $A_j$ is homeomorphic to $S$ and $A_j \subseteq V_0$. Obviously, $j>k$. Thus, $E:=\{ a_{i,k} \}$ is in $\mathcal{D}_i(i+j)$. Let $F:=A_j \setminus \{ a_{j,l} \, ; \, 0 \leq l \leq i+j \}$ and $K := E \times F$. Then $F \in \mathcal{D}_j(i+j)$, $K \in \mathcal{A}^*$ and $K \subseteq U \times A_j \subseteq U \times V_0 \subseteq U \times V \subseteq W$. Moreover, $K$ is homeomorphic to $F$, which is homeomorphic to $A_j$, which is homeomorphic to $S$.

    It is now an immediate consequence of Proposition \ref{ExistenceAndUniqenessOfTheCompactification} that $X \times X$ is homeomorphic to $X$.
\end{proof}

\begin{lemma}\label{UncountableOpenContainsHomeoCopyOfZ(m)}
    Assume $M \subseteq \omega \setminus \{ 0 \}$ is an infinite set and let $S \in \mathscr{S}(M)$. For every uncountable open set $G \subseteq S$, there is a clopen set $H \subseteq S$ contained in $G$ such that $H$ is homeomorphic to $Z(m)$ for some $m \in M$.
\end{lemma}
\begin{proof}
    By the definition of $\mathscr{S}(M)$, we have $S = \alpha_0 \times \dots \times \alpha_k \times Z(m_1) \times \dots \times Z(m_n)$ for some $k,n \in \omega$, $\alpha_0 , \dotsc ,\alpha_k \in \mathcal{O}$ and $m_1, \dotsc ,m_n \in M$. Since $S$ is second-countable and $G$ is uncountable, there is a point in $G$ which does not have any countable neighbourhood. Thus, by the definition of the product topology and by the zero-dimensionality of the spaces $\alpha_0 , \dotsc ,\alpha_k, Z(m_1), \dotsc ,Z(M_n)$, there are clopen sets $U_i \subseteq \alpha_i$, $0 \leq i \leq k$, and clopen sets $V_j \subseteq Z(m_j)$, $1 \leq j \leq n$, such that $W:=U_0 \times \dots \times U_k \times V_1 \times \dots \times V_n$ is uncountable and $W \subseteq G$. Hence, as $\alpha_0 \times \dots \times \alpha_k$ is countable, it follows that $n>0$ and that there is $l \in \{ 1, \dotsc ,n \}$ with $V_l$ uncountable. For every $i \in \{ 0, \dotsc ,k \}$, fix an isolated point $x_i \in U_i$. Since $0 \notin M$, it follows from Lemma \ref{IsolatedPointsAreDenseInZ(alpha)} that there is an isolated point $y_j$ in $V_j$ for every $j \in \{ 1, \dotsc ,n \}$. Let $A_l := V_l$ and $A_j := \{ y_j \}$ for $j \in \{ 1, \dotsc ,n \} \setminus \{ l \}$. Then $H:= \{ x_0 \} \times \dots \times \{ x_k \} \times A_1 \times \dots \times A_n$ is a clopen subset of $S$, it satisfies $H \subseteq W \subseteq G$ and it is homeomorphic to $V_l$. However, $V_l$ is homeomorphic to $Z(m_l)$ by Lemma \ref{UnctblClopenIsHomeo}.
\end{proof}

\begin{proposition}\label{NotHomeomorphicToOpenSubset}
    Let $M \subseteq \omega \setminus \{ 0 \}$ be an infinite set and let $k \in \omega$. Then $Z(k)$ is homeomorphic to an open subset of $X(M)$ if and only if $k \in M$.
\end{proposition}

\begin{proof}
    Again, let us write $X$ instead of $X(M)$, fix a family $\mathcal{A}$ witnessing the defining property of the space $X$ and let $C := X \setminus \bigcup \mathcal{A}$. Obviously, if $k \in M$, then $S:= \{ 0 \} \times Z(k)$ is in $\mathscr{S}(M)$. Hence, fixing a set $A \in \mathcal{A}$ homeomorphic to $S$, we have found a clopen subset of $X$ homeomorphic to $Z(k)$.
    
    Conversely, assume $V$ is an open subset of $X$ homeomorphic to $Z(k)$. If $V \cap C \neq \emptyset$, we just fix arbitrary $m \in M$ and find a set $A \in \mathcal{A}$ with $A \subseteq V$ such that $A$ is homeomorphic to $\{ 0 \} \times Z(m)$, and thus to $Z(m)$. In that case, since $A$ is uncountable and clopen, it follows from Lemma \ref{UnctblClopenIsHomeo} that $A$ is homeomorphic to $Z(k)$, hence $k=m$. Now assume $V \cap C = \emptyset$. Since $V$ is uncountable and $\mathcal{A}$ is countable, there is $A \in \mathcal{A}$ such that $G := V \cap A$ is uncountable. Then, since $A$ is homeomorphic to a member of $\mathscr{S}(M)$, we conclude by Lemma \ref{UncountableOpenContainsHomeoCopyOfZ(m)} that there is a clopen set $H$ contained in $G$ such that $H$ is homeomorphic to $Z(m)$ for some $m \in M$. Then, however, $H$ is an uncountable clopen subset of $V$, and thus it is homeomorphic to $Z(k)$ by Lemma \ref{UnctblClopenIsHomeo}. Thus, $k=m$.
\end{proof}

\begin{proof}[\textbf{Proof of Theorem \ref{MainResult}.}]
    Let $\mathcal{F} := \big\lbrace X(M) \, ; \ M \textup{ is an infinite subset of } \omega \setminus \{ 0 \} \big\rbrace$. Then $|\mathcal{F}|=\mathfrak{c}$ and every member of $\mathcal{F}$ is a compact metrizable zero-dimensional space. By Proposition \ref{X(M)IsHomeoToItsSquare}, $X \times X$ is homeomorphic to $X$ for each $X \in \mathcal{F}$. Finally, it easily follows from Proposition \ref{NotHomeomorphicToOpenSubset} that $\mathcal{F}$ consists of pairwise non-homeomorphic spaces.
\end{proof}

\begin{remark}
The family constructed in the proof of Theorem \ref{MainResult} has the property that for every member $X$ of the family and for every $k \in \omega$, there are (infinitely many) points $x \in X \setminus \mathrm{PK}(X)$ with $\mathrm{CB}(x,X)=k$.  It is thus natural to ask, for each $n \in \omega$, how many compact metrizable zero-dimensional spaces $X$ there are (up to homeomorphism) such that $X$ is homeomorphic to $X \times X$ and $\mathrm{CB}(x,X)\leq n$ for each $x \in X \setminus \mathrm{PK}(X)$. The situation is unclear even for $n=0$, i.e. we have the following question.
\end{remark}

\begin{question}
    How many compact metrizable zero-dimensional spaces $X$ are there (up to homeomorphism) such that $X$ is homeomorphic to $X \times X$ and each point in $X$ is either isolated or belongs to the perfect kernel of $X$?
\end{question}

\section{A correction of \cite[Theorem 3.3]{CharatonikSahan}}

In the final section we present our proof of \cite[Theorem 3.3]{CharatonikSahan}. The original proof, which was largely based on the proof of \cite[Proposition 8.8]{IllanesNadler}, is not entirely correct. Namely, the equality “$\inf \{ d_1 (p,S(X_1)) : p \in \cup_{k=1}^\infty P_k \} = 0$'' in \cite[p.607 line 4]{CharatonikSahan} is not justified. The corresponding part of the proof of \cite[Proposition 8.8]{IllanesNadler} is correct since any infinite set of isolated points in a compact metric space has a zero distance from the set of limit points of the space. However, the set $\cup_{k=1}^\infty P_k$ in the proof of \cite[Theorem 3.3]{CharatonikSahan} is not necessarily a set of isolated points and the set $S(X_1)$, which is the perfect kernel of $X_1$, does not necessarily contain every limit point of $X_1$.

Let us introduce the following notation. Let $\gamma_0 := 1$ and for every nonzero ordinal $\beta$ denote $\gamma_\beta := \omega^\beta+1$. For every $\alpha < \omega_1$ let $\mathscr{S}_\alpha$ be the family of ordinals defined by $\mathscr{S}_\alpha := \{ \gamma_\beta \, ; \ \beta < \alpha \}$.

\begin{lemma}\label{LemmaClopenPartitionsOfCtblOpenSets}
    Let $X$ be a compact metric space and $\alpha< \omega_1$. Assume $G \subseteq X$ is a countable open set such that $\mathrm{CB}(x,X)<\alpha$ for every $x \in G$. Then there is an $X$-clopen partition $\mathcal{A}$ of $G$ such that $\mathcal{A}$ is a null family and every member of $\mathcal{A}$ is homeomorphic to a member of $\mathscr{S}_\alpha$.
\end{lemma}
\begin{proof}
    If $G$ is finite, we just let $\mathcal{A}:=\big\lbrace \{ x \} ; \, x \in G \big\rbrace$. Assume $G$ is infinite (in particular, $G \neq \emptyset$, hence $\alpha>0$) and let $\{ x_n \, ; \ n \in \omega \}$ be an enumeration of $G$. For every $n \in \omega$, fix a compact neighbourhood $K_n$ of $x_n$ contained in $G$ such that $\mathrm{CB}(K_n) = \mathrm{CB}(x_n,X)$. As $G$ is countable, it is zero-dimensional. Hence, since $G$ is open, it easily follows that every point in $G$ admits a neighbourhood base consisting of clopen subsets of $X$. For each $n \in \omega$, let $U_n \subseteq X$ be a clopen set containing $x_n$ such that $U_n \subseteq K_n$ and $\mathrm{diam}(U_n)<2^{-n}$. Let $V_n:=U_n \setminus \bigcup \{ U_k \, ; \, k<n \}$ for every $n \in \omega$ and let $N:=\{ n \in \omega \, ; \ V_n \neq \emptyset \}$. Then $\{ V_n \, ; \, n \in N \}$ is an $X$-clopen partition of $G$ and it is a null family. Moreover, for every $n \in N$, we have $\beta(n) := \mathrm{CB}(V_n) \leq \mathrm{CB}(U_n) \leq \mathrm{CB}(K_n)=\mathrm{CB}(x_n,X)<\alpha$. For each $n \in N$ with $V_n$ infinite, since $V_n$ is homeomorphic to $\omega^{\beta(n)} \cdot k(n)+1$ for some $k(n) \in \omega \setminus \{0\}$, there is a finite $X$-clopen partition $\mathcal{A}_n$ of $V_n$ such that every member of $\mathcal{A}_n$ is homeomorphic to $\omega^{\beta(n)}+1=\gamma_{\beta(n)} \in \mathscr{S}_\alpha$. For each $n \in N$ with $V_n$ finite, let $\mathcal{A}_n:=\big\lbrace \{ x \} ; \, x \in V_n \big\rbrace$ (then every member of $\mathcal{A}_n$ is homeomorphic to $\gamma_0 = \gamma_{\beta(n)} \in \mathscr{S}_\alpha$). Finally, define $\mathcal{A}:=\bigcup \{ \mathcal{A}_n \, ; \, n \in N \}$.
\end{proof}

\begin{proposition}\label{PartitionIntoClopenSetsOfAllRanks}
    Let $(X,d)$ be an uncountable compact metric space and $\alpha$ a countable ordinal. Assume that $\mathrm{CB}(x,X) < \alpha$ for every $x \in X \setminus \mathrm{PK}(X)$ and that
    \[ \mathrm{PK}(X) \subseteq \overline{\{ x \in X \setminus \mathrm{PK}(X) \, ; \ \mathrm{CB}(x,X)=\beta \}} \]
    for every $\beta < \alpha$. Then there is an $X$-clopen partition $\mathcal{A}$ of $X \setminus \mathrm{PK}(X)$ such that:
    \begin{enumerate}[label=(\roman*),font=\textup,noitemsep]
        \item $\mathcal{A}$ is a null family in $(X,d)$;
        \item every member of $\mathcal{A}$ is homeomorphic to a member of $\mathscr{S}_\alpha$;
        \item for every $S \in \mathscr{S}_\alpha$ and every open set $U \subseteq X$ with $U \cap \mathrm{PK}(X) \neq \emptyset$, there is $A \in \mathcal{A}$ such that $A$ is homeomorphic to $S$ and $A \subseteq U$.
    \end{enumerate}
\end{proposition}
\begin{proof}
    If $\alpha =0$, then $X \setminus \mathrm{PK}(X) = \emptyset$ and $\mathscr{S}_\alpha = \emptyset$, hence the choice $\mathcal{A}:=\emptyset$ works. Assume $\alpha > 0$ and let $(\beta_n)_{n \in \omega}$ be a sequence of ordinals less than $\alpha$ such that the set $\{ n \in \omega \, ; \ \beta_n = \beta \}$ is infinite for each $\beta < \alpha$. Fix a sequence $(\varepsilon_n)_{n \in \omega}$ of positive real numbers converging to zero. For every $n \in \omega$, let $F_n \subseteq \mathrm{PK}(X)$ be a finite $\varepsilon_n$-net for $\mathrm{PK}(X)$. We are going to construct (inductively) a sequence $(\mathcal{V}_n)_{n \in \omega}$ of finite families of pairwise disjoint nonempty clopen subsets of $X$ such that the following conditions hold for every $n \in \omega$.
    \begin{enumerate}[label=(\arabic*),font=\textup,noitemsep]
        \item $\forall \, k < n : \big( \bigcup \mathcal{V}_k \big) \cap \big( \bigcup \mathcal{V}_n \big) = \emptyset$;
        \item $\forall \, V \in \mathcal{V}_n : V \subseteq X \setminus \mathrm{PK}(X)$;
        \item $\forall \, V \in \mathcal{V}_n : \mathrm{CB}(V)=\beta_n$;
        \item $\forall \, V \in \mathcal{V}_n : \mathrm{diam}_{d} (V) < \varepsilon_n$;
        \item $\forall \, V \in \mathcal{V}_n : \mathrm{dist}_d (\mathrm{PK}(X),V) < \varepsilon_n$;
        \item $\forall \, z \in F_n \ \exists \, V \in \mathcal{V}_n : \mathrm{dist}_d (z,V) < \varepsilon_n$.
    \end{enumerate}
    Before we start the construction, for every $x \in X \setminus \mathrm{PK}(X)$, let $K_x$ be a compact neighbourhood of $x$ such that $K_x \subseteq X \setminus \mathrm{PK}(X)$ and $\mathrm{CB}(K_x)=\mathrm{CB}(x,X)$. Now, let us construct $\mathcal{V}_0$. For every $z \in F_0$, since $\mathrm{PK}(X) \subseteq \overline{\{ x \in X \setminus \mathrm{PK}(X) \, ; \ \mathrm{CB}(x,X)=\beta_0 \}}$, there is $x_0(z) \in X \setminus \mathrm{PK}(X)$ such that $\mathrm{CB} \big( x_0(z),X \big)=\beta_0$ and $d \big( z, x_0(z) \big) < \varepsilon_0$. Clearly, we can assume that $x_0(z) \neq x_0 (w)$ for every two distinct points $z,w \in F_0$. Since $X \setminus \mathrm{PK}(X)$ is countable (and thus zero-dimensional) and open, every point in $X \setminus \mathrm{PK}(X)$ admits a neighbourhood base consisting of clopen subsets of $X$. Thus, for every $z \in F_0$, there is a clopen set $V_0(z) \subseteq X$ with $x_0(z) \in V_0(z) \subseteq K_{x_0(z)}$ and $\mathrm{diam}_d \big( V_0(z) \big) < \varepsilon_0$. Again, we can assume that $V_0(z) \cap V_0(w) = \emptyset$ for every two distinct points $z,w \in F_0$. For every $z \in F_0$, we have $\mathrm{dist}_d \big( \mathrm{PK}(X), V_0(z) \big) \leq \mathrm{dist}_d \big( z, V_0(z) \big) \leq d \big( z, x_0(z) \big) < \varepsilon_0$ and
    \[ \beta_0 = \mathrm{CB} \big( x_0(z), X \big) \leq \mathrm{CB} \big( V_0(z)  \big) \leq \mathrm{CB} \big( K_{x_0(z)} \big) = \beta_0. \]
    Hence, letting $\mathcal{V}_0 := \big\lbrace V_0(z) \, ; \, z \in F_0 \big\rbrace$, the initial step is done. In the same fashion, given $n \in \omega$ and assuming the families $\mathcal{V}_0, \dotsc , \mathcal{V}_n$ have already been constructed, it is possible to construct $\mathcal{V}_{n+1}$. Note that it is easy to make sure (1) is satisfied as $\bigcup \big\lbrace \bigcup \mathcal{V}_i \, ; \ i \leq n \big\rbrace$ is a compact set disjoint from $\mathrm{PK}(X)$.

    Having finished the inductive construction, let $\mathcal{V} := \bigcup \{ \mathcal{V}_n \, ; \ n \in \omega \}$. Then $\mathcal{V}$ is a null family of pairwise disjoint nonempty clopen subsets of $X$ contained in $X \setminus \mathrm{PK}(X)$. For every $V \in \mathcal{V}$, there is (by the same argument as in the proof of Lemma \ref{LemmaClopenPartitionsOfCtblOpenSets}) a finite $X$-clopen partition $\mathcal{A}_V$ of $V$ such that every member of $\mathcal{A}_V$ is homeomorphic to $\gamma_{\mathrm{CB}(V)} \in \mathscr{S}_\alpha$. It easily follows from (4) and (5) that the set $G:= \big( X \setminus \mathrm{PK}(X) \big) \setminus \bigcup \mathcal{V}$ is open in $X$. Thus, by Lemma \ref{LemmaClopenPartitionsOfCtblOpenSets}, there is an $X$-clopen partition $\mathcal{A}_G$ of $G$ such that $\mathcal{A}_G$ is a null family and every member of $\mathcal{A}_G$ is homeomorphic to a member of $\mathscr{S}_\alpha$. Let $\mathcal{A}:= \mathcal{A}_G \cup \bigcup \{ \mathcal{A}_V \, ; \, V \in \mathcal{V} \}$. Clearly, $\mathcal{A}$ is an $X$-clopen partition of $X \setminus \mathrm{PK}(X)$ and it satisfies (i) and (ii). To verify assertion (iii), let $S \in \mathscr{S}_\alpha$ be given and let $U \subseteq X$ be an open set intersecting $\mathrm{PK}(X)$. Obviously, $S=\gamma_\beta$ for some $\beta < \alpha$. Fix a point $p \in U \cap \mathrm{PK}(X)$ and let $r$ be a positive real number such that the open $r$-ball centered at $p$ is contained in $U$. There is $m \in \omega$ such that $\varepsilon_n < r/3$ for each $n \in \omega$ with $n\geq m$. Since the set $\{ n \in \omega \, ; \ \beta_n = \beta \}$ is infinite, there is $n \in \omega$ with $n\geq m$ such that $\beta_n = \beta$. As $F_n$ is an $\varepsilon_n$-net for $\mathrm{PK}(X)$, there is $z \in F_n$ such that $d(p,z)<\varepsilon_n$. By (6), there is $V \in \mathcal{V}_n$ with $\mathrm{dist}_d (z,V) < \varepsilon_n$. Then $\mathrm{dist}_d (p,V) < 2\varepsilon_n < 2r/3$ and thus, since $\mathrm{diam}_{d} (V) < \varepsilon_n <r/3$ by (4), the set $V$ is contained in the open $r$-ball centered at $p$. This shows that $V \subseteq U$. Moreover, $\mathrm{CB}(V)=\beta_n=\beta$ by (3). Hence, taking any $A \in \mathcal{A}_V$, we have found a member of $\mathcal{A}$ contained in $U$ and homeomorphic to $\gamma_\beta = S$.
\end{proof}

Combining Proposition \ref{ExtendingHomeos} with Proposition \ref{PartitionIntoClopenSetsOfAllRanks}, we immediately obtain the following restatement of \cite[Theorem 3.3]{CharatonikSahan}.

\begin{corollary}
    Let $\alpha$ be a countable ordinal and let $X_1$, $X_2$ be uncountable metrizable compact spaces such that the following conditions are satisfied for each $i \in \{ 1,2 \}$.
    \begin{itemize}
        \item $\mathrm{CB} (x,X_i) < \alpha$ for every $x \in X_i \setminus \mathrm{PK} (X_i)$;
        \item $\mathrm{PK} (X_i) \subseteq \overline{\big\lbrace x \in X_i \setminus \mathrm{PK} (X_i) \, ; \ \mathrm{CB} (x,X_i) = \beta \big\rbrace}$ for every $\beta < \alpha$.
    \end{itemize}
    Then every homeomorphism $h \colon \mathrm{PK} (X_1) \to \mathrm{PK} (X_2)$ extends to a homeomorphism $\overline{h} \colon X_1 \to X_2$.
\end{corollary}

\bibliographystyle{alpha}
\bibliography{citace}
\end{document}